\documentclass[runningheads]{llncs}

\usepackage{graphicx}
\usepackage{lineno,hyperref}
\usepackage{booktabs}
\usepackage{amsmath,amssymb,amsfonts}
\usepackage{amssymb}
\usepackage{eqnarray} 
\usepackage{mathalfa}  
\usepackage{algorithm}
\usepackage{algpseudocode}
\usepackage{psfrag}
\usepackage{stfloats}
\usepackage{dsfont}
\usepackage{xcolor}
\usepackage{bbm}
\usepackage{float}
\usepackage{comment}
\usepackage[prependcaption,colorinlistoftodos]{todonotes}
\usepackage{hyperref}
\usepackage[per-mode=repeated-symbol]{siunitx}
\usepackage{subfigure}

\DeclareMathOperator{\argmin}{\operatorname{arg min}}

\newcommand{\sF}{\mathcal{F}}

\newcommand{\sX}{\mathcal{X}}
\newcommand{\sU}{\mathcal{U}}

\newcommand{\sW}{\mathcal{W}}

\newcommand{\R}{\mathbb{R}}
\newcommand{\E}{\mathbb{E}}

\newcommand{\bv}[1]{\mathbf{#1}}

\newcommand{\subscr}[2]{#1_{\textup{#2}}}

\newcommand{\np}{\subscr{n}{$\pi$}}

\newcommand{\marginalxk}[1]{p_{{#1}}\left(\bv{x}_{{#1}}\right)}
\newcommand{\refmarginalxk}[1]{q_{{#1}}\left(\bv{x}_{{#1}}\right)}

\newcommand{\jointxuk}[2]{p_{{#1}}\left(\bv{x}_{{#1}},\bv{u}_{{#1}}{\mid} \bv{x}_{{#2}} \right)}
\newcommand{\refjointxuk}[2]{q_{{#1}}\left(\bv{x}_{{#1}},\bv{u}_{{#1}}{\mid} \bv{x}_{{#2}} \right)}

\newcommand{\plantk}[2]{p_{\bv{x},{#1}} \left(\bv{x}_{{#1}}{\mid} \bv{x}_{{#2}}, \bv{u}_{{#1}} \right)}
\newcommand{\refplantk}[2]{q_{\bv{x},{#1}} \left(\bv{x}_{{#1}}{\mid} \bv{x}_{{#2}}, \bv{u}_{{#1}} \right)}

\newcommand{\policyk}[2]{p_{\bv{u},{#1}}\left(\bv{u}_{{#1}}{\mid} \bv{x}_{{#2}} \right)}
\newcommand{\refpolicyk}[2]{q_{\bv{u},{#1}}\left(\bv{u}_{{#1}}{\mid} \bv{x}_{{#2}} \right)}
\newcommand{\primitivek}[3]{\pi_{{#1}}^{({#3})}\left(\bv{u}_{{#1}}{\mid} \bv{x}_{{#2}} \right)}

\newcommand{\marginalx}[1]{p\left(\bv{x}_{{#1}}\right)}

\newcommand{\jointxu}[2]{p\left(\bv{x}_{{#1}},\bv{u}_{{#1}}{\mid} \bv{x}_{{#2}} \right)}
\newcommand{\refjointxu}[2]{q\left(\bv{x}_{{#1}},\bv{u}_{{#1}}{\mid} \bv{x}_{{#2}} \right)}

\newcommand{\plant}[2]{p_{\bv{x}}\left(\bv{x}_{{#1}}{\mid} \bv{x}_{{#2}}, \bv{u}_{{#1}} \right)}
\newcommand{\refplant}[2]{q_{\bv{x}} \left(\bv{x}_{{#1}}{\mid} \bv{x}_{{#2}}, \bv{u}_{{#1}} \right)}

\newcommand{\policy}[2]{p_{\bv{u}}\left(\bv{u}_{{#1}}{\mid} \bv{x}_{{#2}} \right)}
\newcommand{\refpolicy}[2]{q_{\bv{u}}\left(\bv{u}_{{#1}}{\mid} \bv{x}_{{#2}} \right)}
\newcommand{\primitive}[3]{\pi^{({#3})}\left(\bv{u}_{{#1}}{\mid} \bv{x}_{{#2}} \right)}

\newcommand{\weight}[2]{\bv{w}_{#1}^{(#2)}}
\newcommand{\weights}[1]{\bv{w}_{#1}}
\newcommand{\optimalweights}[1]{\bv{w}^{\star}_{#1}}

\newcommand{\costx}[1]{c_{#1}^{\textup{x}}\left(\bv{X}_{#1}\right)}
\newcommand{\costu}[1]{c_{#1}^{\textup{u}}\left(\bv{U}_{#1}\right)}
\newcommand{\costl}[2]{l_{#1}^{\star}\left(\bv{X}_{#2}\right)}

\newcommand{\optimalpolicyk}[2]{p^{\star}_{\bv{u},{#1}} \left(\bv{u}_{{#1}}{\mid} \bv{x}_{{#2}} \right)}

\newcommand{\optimalpolicy}[2]{p^{\star}_{\bv{u}} \left(\bv{u}_{{#1}}{\mid} \bv{x}_{{#2}} \right)}
\newcommand{\optimalweight}[2]{\bv{w}^{(#2), \star}_{#1}}

\newcommand{\KL}{\textup{KL}}
\newcommand{\DKL}[2]{{D}_{\KL}\left(#1{\mid \mid} #2 \right)}

\newtheorem{Remark}{Remark}

\newtheorem{Property}{Property}

\newtheorem{Assumption}{Assumption}

\newtheorem{Problem}{Problem}

\newcommand{\Model}{Free-Gate}

\newcommand{\ds}{\displaystyle}
\usepackage[T1]{fontenc}
\usepackage{graphicx}
\usepackage{times}

\begin{document}
\title{\Model: Planning, Control And Policy Composition via Free Energy Gating}

\titlerunning{Planning, Control And Policy Composition via Free Energy Gating}
%

\author{Francesca Rossi\inst{1}\orcidID{0009-0005-5663-260X} \and
{\'E}miland Garrab{\'e}\inst{2}\orcidID{0000-0003-4775-0258} \and
Giovanni Russo$^*$\inst{1,3}\orcidID{0000-0001-5001-3027}}

\authorrunning{F. Rossi et al.}

\institute{Scuola Superiore Meridionale, Naples, Italy \email{f.rossi@ssmeridionale.it} 
\and
Sorbonne University, Paris, France \email{garrabe@isir.upmc.fr} 
\and
University of Salerno, Italy \email{giovarusso@unisa.it}}

\maketitle              
\begin{abstract}
We consider the problem of optimally composing a set of primitives to tackle planning and control tasks. To address this problem, we introduce a free energy computational model for planning and control via policy composition: \Model. Within \Model, control primitives are combined via a gating mechanism that minimizes variational free energy. This composition problem is formulated as a finite-horizon optimal control problem, which we prove remains convex even when the cost is not convex in states/actions and the environment is nonlinear, stochastic and non-stationary. We develop an algorithm that computes the optimal primitives composition and demonstrate its effectiveness via in-silico and hardware experiments on an application involving robot navigation in an environment with obstacles. The experiments highlight that \Model~enables the robot to navigate to the destination despite only having available simple motor primitives that, individually, could not fulfill the task.

\keywords{Free Energy Minimization \and Policy Composition \and Autonomous Systems \and Decision Making \and  Optimal Control }
\end{abstract}

\section{Introduction}
When we are given a new task, we can re-use and compose what we know to achieve our goal. This ability to compose knowledge is widely believed to be a key ingredient in differentiating natural and machine intelligence \cite{lake2017building,russin2024,Vijayaraghavan2025}. For state-of-the-art autonomous agents, designed by, e.g., leveraging deep neural networks, this desirable capability is promoted by endowing the network with inductive biases. In turn, this is achieved via, e.g., specialized architectures, meta-learning, or large-scale pre-training \cite{russin2024}. However, there is little evidence that machines designed on this paradigm can achieve knowledge composition and a key challenge transversal to learning and control is to equip autonomous agents with mechanisms mimicking this ability.

Motivated by this, we introduce \Model: a computational model for planning and control through policy composition. In \Model, simple policies (i.e., primitives) are composed by a gating mechanism. This mechanism is grounded in the minimization of variational free energy, a unifying account across information theory, statistical learning and active inference \cite{OS:18,BZ_AM_AB_AD:08,friston2010free}. Essentially, \Model~enables agents to compute an optimal policy from a set of primitives that are dynamically composed via a variational free energy minimization process. In turn, this composition problem is formalized as a finite-horizon optimal control problem. \Model~returns the optimal solution to this problem: the (possibly, non-stationary) optimal policy for the agent built by linearly combining the primitives. After characterizing convexity of the problem we give an algorithm that provably returns the optimal solution. The effectiveness of the algorithm is demonstrated on a robot navigation task in an environment with obstacles. Remarkably, our in-silico and hardware experiments show that \Model~enables  the robot to navigate to the goal position -- and avoid obstacles -- despite only having available primitives that, individually, could not complete the task. Next, we survey policy composition architectures related to \Model~(see also Section \ref{sec:background}) and refer interested readers to \cite{lake2017building,russin2024} for broader discussions on compositionality in natural and machine intelligence from the viewpoint of psychology, cognition and philosophy.

\paragraph{\textbf{Related Works.}} The design of modular architectures enabling agents to fulfill a given task from a set of primitives is a central topic for learning and control, with robotics often used as a design, test and validation domain \cite{prescott2023understanding,Vijayaraghavan2025}. In this context, two popular models are MOSAIC (Modular Selection and Identification for Control) and the mixture of experts, see, e.g., \cite{mosaic,haruno2001mosaic,jacobs1991adaptive,van2022embodied}. The MOSAIC architecture builds control inputs via a weighted linear combination of motor {\em command} signals. The weights for the commands are obtained from {\em responsibility signals} computed via a stationary soft-max selection rule. Building on this,  \cite{doya2002multiple,samejima2006multiple} complement the soft-max rule for weights selection with temporal difference learning aimed at maximizing a stationary reward.  Despite the experimental successes on nonlinear control tasks reported in these works, to the best of our knowledge there is no normative framework explaining why the soft-max rule should be used to compose the motor signals. The mixture of experts model was originally introduced in the context of machine learning \cite{jacobs1991adaptive} and has become a cornerstone for Large Language Model architectures \cite{cai2025survey}. Following this approach, decisions are computed by learning a stationary policy (in infinite-horizon settings) to linearly combine the output of a pool of experts that also set the performance baseline. We refer to, e.g., \cite{NCB_GL:06,masoudnia2014mixture} for more details and to \cite{cacciatore1993mixtures,gomi1993recognition,willi2024mixture} for applications to control tasks where a stationary policy is again learned for infinite-horizon control problems. 
A related approach to combine experts is via architectures based on the so-called {Product of Experts \cite{hinton2002training} model where experts (probabilistic models) are combined by properly multiplying their probability distributions. While originally developed for inference, this model has been recently extended \cite{hansel2023hierarchical} to learn stationary policies for control tasks. Non-stationary policies (in finite-horizon settings) can be instead computed with the architecture from \cite{garrabe2023optimal}. However, the results are obtained under the simplifying assumption that the agent can directly specify state transitions rather than policies.
Policy composition is also related to the rich literature on hierarchical Reinforcement Learning (RL) architectures~\cite{pateria2021hierarchical}. The idea behind these architectures is to break down long time-horizon (and complex)  tasks into a sequence of shorter (and simpler) subtasks. Namely, a {\em high-level} policy is learned to orchestrate the subtasks and each of the subtasks is in turn learned via its own RL algorithm, thus yielding a \textit{temporal abstraction}~\cite{pateria2021hierarchical,vezhnevets2017feudal}. The hierarchical decomposition can shorten the task learning time, thus improving learning efficiency~\cite{hutsebaut2022hierarchical}. Remarkably, in~\cite{pezzulo2018hierarchical} this approach is applied to active inference, with behavior emerging from the minimization of expected free energy at each level of the hierarchy. {In the context of active inference we also recall~\cite{Paul2024}. Here, hierarchical architectures are proposed in which multiple generative models are selectively engaged based on contextual cues, thus implementing a mixture of experts scheme guided by epistemic and instrumental value. These models enable efficient planning by composing both perception and action strategies in a modular way. In \cite{Paul2024b} connections are drawn between active inference and KL control and this enables efficient computation (see also \cite{EG_GR:22} for a survey across learning and control). In developmental robotics, compositionality emerges through interactive learning processes that couple sensory and motor experiences with symbolic structure~\cite{tani1999learning,Vijayaraghavan2025}. In these works world models are constructed by composing simple generative components, supporting abstraction and behavior generalization.} In the context of computational neuroscience, an architecture for knowledge composition is also postulated within the Thousand Brains Theory (TBT). TBT argues \cite{hawkins2021thousand,Mountcastle1978} that our ability to efficiently learn new concepts is hosted in our neocortex, a highly parallel system having cortical columns as functional units. According to TBT, our perception of the world arises from composing the outputs of these parallel units.

\paragraph{\textbf{Contributions.}} We present \Model, a computational model for planning and control enabling agents to optimally combine control primitives in order to tackle sequential tasks. In \Model, primitives are optimally composed via a gating mechanism that minimizes variational free energy. This leads to formalize a finite-horizon optimal control problem that is convex even when the cost is not convex in states/actions and the environment is nonlinear, non-stationary and stochastic. To the best of our knowledge, this is the first free energy model featuring an underlying optimization problem for optimal policy computation that is convex in this setting and that, additionally: (i) does not require that the agent {directly} specifies its state transitions; (ii) does not assume stationary/discounted settings. Our key technical contributions are as follows: 
\begin{itemize}
    \item our model provides a normative mechanism for primitives composition. This mechanism links the decision rules determining how primitives are composed to the solution of an underlying constrained optimal control problem having the variational free energy as cost and primitives' weights as decision variables;
    \item we  characterize the underlying optimization problem. This can be tackled via recursive arguments and, at each time step, the optimal weights for primitives composition can be found by solving a problem that is convex even when the agent cost is not in states/actions. We do not make specific assumptions on the environment dynamics besides the standard Markov assumption;
    \item the results are turned into an algorithm and validated both in-silico, via simulations, and on a hardware test-bed, using the Robotarium \cite{robotarium}. The experiments show the effectiveness of \Model~on a rover navigation task, where the agent needs to combine simple motor primitives to reach a final destination while avoiding obstacles.  Finally, the code implementing \Model~is fully documented  and made available at \url{https://github.com/freegate-iwai/Free-Gate.git}.
\end{itemize}
The paper is organized as follows. After introducing (Section \ref{sec:background}) the background, we present  \Model~and recast the primitives composition problem as an optimal control problem (Section \ref{sec:architecture}). Then, we present an algorithm (Section \ref{sec:properties}) to tackle the problem, characterizing: (i)  optimality; (ii) convexity of the underlying optimization. \Model~is validated in Section \ref{sec:validation} and concluding remarks are given in Section \ref{sec:conclusions}.

\section{Background}\label{sec:background}
 Sets are denoted in calligraphic, e.g., $\mathcal{X}$, and vectors in {\textbf{bold}}, e.g., $\mathbf{x}$. A random variable is denoted by $\bv{V}$ and its realization by $\bv{v}$. For continuous (or discrete) random variables we denote the \textit{probability density function}, pdf, (or \textit{probability mass function}, pmf) of $\mathbf{V}$ by $p(\mathbf{v})$. The joint pmf of $\mathbf{V}_1$ and $\mathbf{V}_2$ is denoted by $p(\mathbf{v}_1,\mathbf{v}_2)$ and the conditional pmf of $\mathbf{V}_1$ with respect to $\mathbf{V}_2$ is $p\left( \mathbf{v}_1\mid \mathbf{v}_2 \right)$. Countable sets are denoted by $\lbrace w_k \rbrace_{k_1:k_n}$, where $w_k$ is a generic element, $k_1$ and $k_n$ are respectively the indices of the first and last element, and $k_1:k_n$ is the set of consecutive integers between (including) $k_1$ and $k_n$. The support of a function $p(\cdot)$ is $\mathcal{S}(p)$. The expectation of a function $\mathbf{h}(\cdot)$ of a continuous random variable $\mathbf{V}$ with pdf $p(\mathbf{v})$ is $\E_{{p}}[\mathbf{h}(\mathbf{V})]:=\int \mathbf{h}(\mathbf{v}) p(\mathbf{v}) d\bv{v}$. We utilize the Kullback-Leibler (KL) divergence \cite{SK_RL:51}: $ \DKL{p}{q}:= \int p(\bv{v}) \ln\left( \frac{p(\bv{v})}{q(\bv{v})}\right)d\bv{v}$, that is a measure of discrepancy between two pdfs $p(\mathbf{v})$ and $q(\mathbf{v})$. The $\DKL{p}{q}$ is finite only if $\mathcal{S}(p)\subset \mathcal{S}(q)$, i.e., $p(\bv{v})$ is absolutely continuous with respect to $q(\bv{v})$, see \cite[Chapter $8$]{TC_JT:06}. We use the standard convention $0\log(0) = 0$. For discrete variables the integral is replaced with the sum in the above expressions. {Finally, we let
$\Delta^{n} := \left\{ \bv{w} \in \mathbb{R}^{n} \;\middle|\; w_{i} \geq 0 \ \forall i \in 1:n,\; \sum_{i=1}^{n} w_i = 1 \right\}$ be the probability simplex and recall that the affine dimension of $\Delta^{n} \subset \mathbb{R}^{n}$ is $n - 1$.}

\paragraph{\textbf{Free Energy Minimization.}} The minimization of free energy plays a central role for autonomous decision making across optimization, inference, risk-aware control and learning \cite{PJ_HH_VR:24,TS:14,MR_etal:21,kouw2023information,EG_GR:22}. A wide range of methods and frameworks, including RL schemes \cite{BZ_AM_AB_AD:08,haarnoja2017reinforcement,Garrabe2025} based on maximum entropy, Bayesian inference \cite{MJ_ZG_TJ_LS:99}, and entropy-regularized optimal transport \cite{GC_VD_GP_BS:17}, rely on optimizing free energy functionals \cite{STJ_OS:21,OS:18}, eventually with constraints \cite{DG_GR:22}. In the context of active inference -- which unifies perception, action, and learning -- the variational free energy is minimized to align the agent internal model with sensory inputs \cite{friston2010free,friston2017active}. For control, this framework optimizes policies under uncertainty using actions as decision variables while maintaining consistency with internal models \cite{millidge2020relationship,pezzulo2015active,Paul2024}. Beyond artificial agents, this minimization principle is also central to the broader quest for a unifying theory of self-organization \cite{CH_etal:24}, brain function and sentient behaviors \cite{KF:09,KF_LDC_NS_CH_KU_GP_AG_TP:23,Parr2022}. Formally, see, e.g.,~\cite{STJ_OS:21}, the free energy is a functional of the form:

\begin{equation*}
\sF(p(\bv{z})) = \DKL{p(\bv{z})}{q(\bv{z})} + \E_{p(\bv{z})} \left[ l(\bv{z}) \right],
\end{equation*}
where: (i) the first term is the statistical complexity of $p(\bv{z})$ with respect to a given $q(\bv{z})$. Minimizing this term amounts at minimizing the discrepancy (in the KL divergence sense) between $p(\bv{z})$ and $q(\bv{z})$; (ii) the second term is the agent's expected loss. Following, e.g.,~\cite{STJ_OS:21}, when (as in active inference) $l(\bv{z})$ is a log-loss of the form $-\log(\tilde{p}(\bv{z}))$ the free energy is known as the evidence lower bound (ELBO) commonly used in variational inference \cite[Chapter 10]{KM:23}.

\begin{figure}[tb]
    \centering
    \includegraphics[width=0.7\columnwidth]{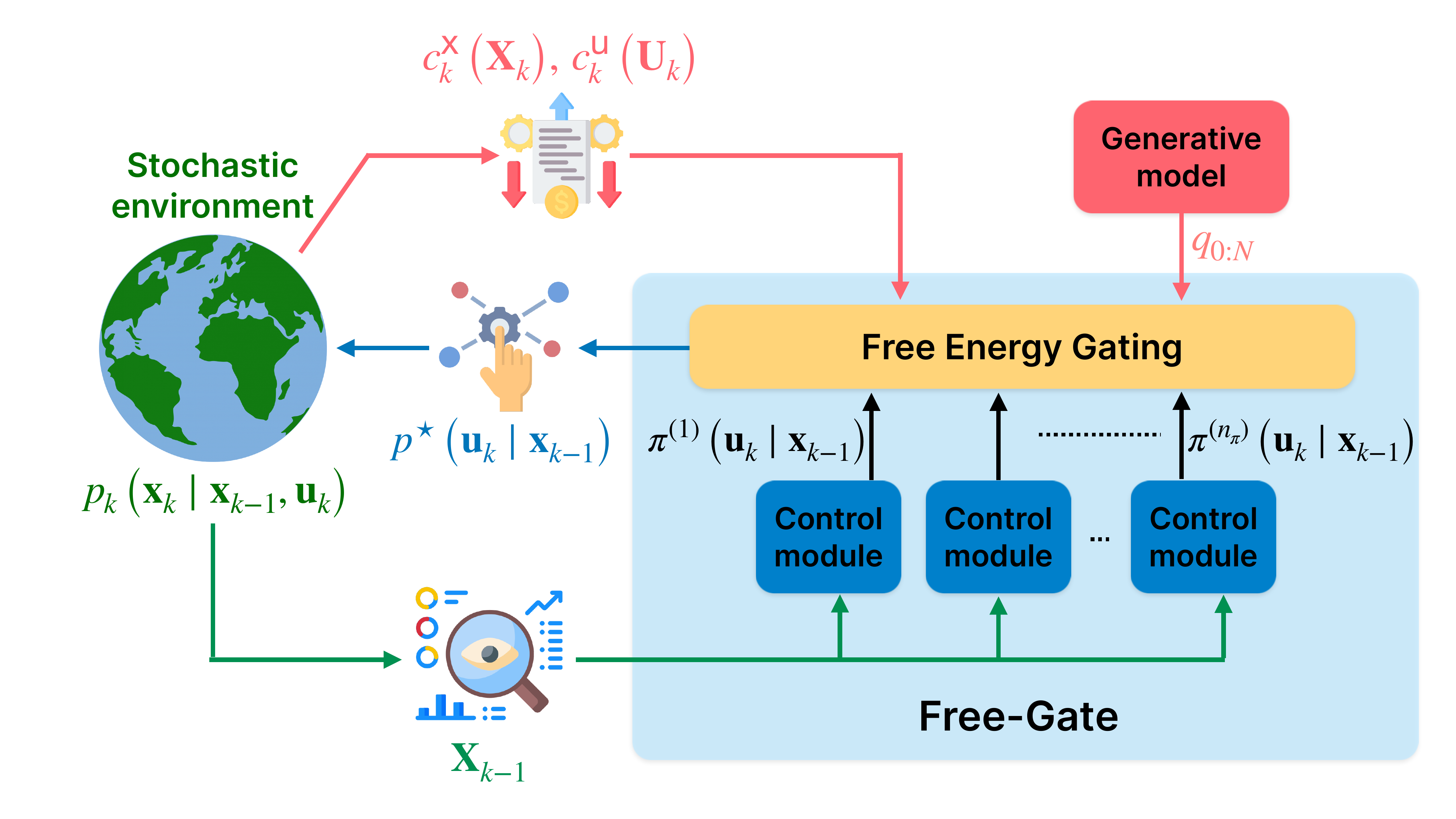}
    \caption{In \Model~control primitives are dynamically combined by a free energy gating mechanism. This mechanism, given the state, action/state costs and the generative model, computes the optimal weights for the primitives composition by minimizing variational free energy. Consequently, the \Model~possibly non-stationary policy, $\optimalpolicyk{k}{k-1}$, is a linear combination of these primitives with the optimal weights determined via free energy gating.}\label{fig:architecture_b}
\end{figure}
\noindent

\section{The \Model~Model}

We consider the set-up of Fig. \ref{fig:architecture_b} where the agent determines an action $\bv{u}_k\in\sU\subseteq\R^{n_u}$ based on the state $\bv{x}_{k-1}\in\sX\subseteq\R^{n_x}$. The time indexing is chosen so that the environment/system transitions from $\bv{x}_{k-1}$ to $\bv{x}_k$ when $\bv{u}_k$ is applied. The possibly stochastic, nonlinear and non-stationary Markovian environment is described, at each $k$, by $\plantk{k}{k-1}$. The agent determines the action by sampling from the policy $\policyk{k}{k-1}$ so that the closed-loop behavior over $k\in0:N$ is described by $p_{0:N} := \marginalxk{0} \ds\prod_{k=1}^N\jointxuk{k}{k-1}$, with $\marginalxk{0}$ being an initial prior capturing initial conditions and $\jointxuk{k}{k-1} := \plantk{k}{k-1}\policyk{k}{k-1}$. Also, $c_{k}^{\textup{x}}(\cdot):\sX\rightarrow \R$ and $c_{k}^{\textup{u}}(\cdot):\sU\rightarrow \R$ are the state and action costs, respectively.

\subsection{\Model~Description}\label{sec:architecture}

As schematically illustrated in Fig. \ref{fig:architecture_b}, given state/action costs and a generative model
$$
q_{0:N}{:=} \refmarginalxk{0}\prod_{k=1}^N\refjointxuk{k}{k-1}{=} \refmarginalxk{0}\prod_{k=1}^N\refplantk{k}{k-1}\refpolicyk{k}{k-1},
$$
the aim of~\Model~is to compute $\{\optimalpolicyk{k}{k-1}\}_{1:N}$, an optimal sequence of policies that minimizes the variational free energy functional 
$$
\DKL{p_{0:N}}{q_{0:N}} + \sum_{k=1}^N\E_{\marginalxk{k-1}}\left[\E_{\jointxuk{k}{k-1}}\left[\costx{k} + \costu{k}\right]\right],
$$
where $\marginalxk{k-1}$ is the marginal probability of the state at $k-1$. The functional naturally arises in the context of control and learning (see, e.g., \cite{AS_HJ_KF_GR:25,Garrabe2025}) and, when the cost is the negative log-likelihood of states, it yields the variational free energy minimized in active inference based on the expected free energy~\cite{AS_HJ_KF_GR:25,Parr2019}. At each $k$, $\optimalpolicyk{k}{k-1}$ is a linear combination of a set of $\np$ primitives, denoted by $\primitivek{k}{k-1}{i}$, $i \in \mathcal{P} := 1:\np$ and have support $\mathcal{S}(\pi)$. In the architecture: (i) each {\em 
{control module}} provides a primitive. The primitive is black-box, in the sense that the {control module} does not {\em disclose} what underlying control algorithm it is implementing. The terminology {\em primitive} is inspired by the fact that the {modules} provide policies encoding simple behaviors (e.g., basic motor commands as in our experiments); (ii) {\em free energy gating} mechanism, which given the primitives returns the weights, $\optimalweight{k}{i}$, to optimally combine them. Thus, the optimal policy is $\optimalpolicyk{k}{k-1} = \ds \sum_{i\in \mathcal{P}}\optimalweight{k}{i}\primitivek{k}{k-1}{i}$. \\

\noindent{\bf Free Energy Gating.} In \Model, the gating mechanism aims at finding the optimal primitives' weights combination by minimizing variational free energy over the time horizon $1:N$. Namely, the problem of finding the optimal weights combination is recast via the following finite-horizon optimal control problem (to streamline notation, from now on we omit time steps in the probabilities subscripts).
\begin{Problem}\label{prob:finite-horizon}
Given {the $\{\plant{k}{k-1}\}_{1:N}$,} the generative model $q_{0:N}$, and the set of primitives $\primitive{k}{k-1}{i}$, $i\in \mathcal{P}$, find the sequence of weights $\{\optimalweights{k}\}_{1:N}$ such that:
\begin{equation}\label{eqn:finite-horizon}
    \begin{aligned}
    {\left\{\optimalweights{k}\right\}_{1:N}} {\in} \underset{\left\{ \weights{k} \right\}_{1:N}}{\argmin} 
    &\DKL{p_{0:N}}{q_{0:N}}{+}\!\! \sum_{k=1}^N\E_{\marginalx{k-1}}{\left[\E_{\jointxu{k}{k-1}}{\left[\costx{k}\!{+}\costu{k}\right]}\right]}\\
    \mbox{s.t. } & \policy{k}{k-1} = \sum_{i\in \mathcal{P}}\weight{k}{i}\primitive{k}{k-1}{i}, \ \ \ \forall k\in 1:N \\
    & \weights{k} \in \Delta^{\np}, \ \ \ \forall  k \in 1:N,
    \end{aligned}
\end{equation}
where $\weights{k}$ is the stack of $\weight{k}{i}$'s and $\marginalx{k-1}$ is the state marginal probability  at $k-1$.
\end{Problem}
In the formulation, the costs are left general and can be used to encode the agent task. For example, in our robot navigation application, the costs promote obstacle avoidance. In Problem \ref{prob:finite-horizon}, the constraints formalize the fact that the agent policy is a combination of the primitives. The cost functional, which is also considered in~\cite{Garrabe2025}, becomes the surprise upper bound typical in active inference when the cost terms inside the expectation are chosen as negative log-likelihoods~\cite{AS_HJ_KF_GR:25}. We also refer to, e.g.,~\cite{Paul2024b} and~\cite{Kouw2024} for an application to planning, and control from demonstration, where the generative model encodes a desired behavior extracted from example data~\cite{DG_GR:22}. The decision variables are the weights of the linear combination and therefore, at each $k$, the \Model~policy is $\optimalpolicy{k}{k-1}=\sum_{i\in \mathcal{P}}\optimalweight{k}{i}\primitive{k}{k-1}{i}$. In Problem \ref{prob:finite-horizon}, $N = 1$ yields biomimetic policies \cite{vincent2006biomimetics} and these often arise as reflexes in bio-inspired motor control \cite{pezzulo2015active}. Next, we let $\bv{w}$ be the stack of all $\bv{w}_k$'s, {$k {\in} 1:N$}; $F_{0:N}(\bv{w})$ denotes the cost of Problem \ref{prob:finite-horizon}. We make the following standing assumptions:
\begin{Assumption}\label{assum_primitives}
    The primitives are bounded.
\end{Assumption}
\begin{Assumption}\label{asn:bounded_cost}
$\sW := \big\{\bv{w} \colon \bv{w} \ \mbox{is feasible for Problem \ref{prob:finite-horizon} and} \ F_{0:N}(\bv{w}) < +\infty  \big\}$ is non-empty.
\end{Assumption}
\begin{Assumption}\label{assum_support}
    Each primitive has full support over the action space, i.e., $\mathcal{S}(\pi)=\mathcal{U}$.
\end{Assumption}
\begin{Remark}\label{rem:assumptions}
Assumption \ref{assum_primitives} is rather mild and always satisfied for, e.g., discrete variables. Assumption \ref{asn:bounded_cost} is standard in the context of optimal transport, see, e.g., \cite[Theorem 1.10]{MN_20}. This assumption, guaranteeing that the optimal cost is finite, is satisfied when: (i) the state/action cost functions are lower-bounded; (ii) $\plant{k}{k-1}$ is absolutely continuous with respect to $\refplant{k}{k-1}$; (iii) the primitives are absolutely continuous with respect to $\refpolicy{k}{k-1}$. {Assumption \ref{assum_support} is used to guarantee that the gradient of the cost in~\eqref{eqn:finite-horizon} is well defined. This assumption, implying that the policy is not deterministic, can also be found in the context of policy gradient methods for RL (see, e.g.,~\cite[Chapter~13]{sutton1998reinforcement}).}
\end{Remark}

\begin{algorithm}[tb]
    \caption{Free Energy Gating}\label{alg_s2}
    \begin{algorithmic}[1]
\State \textbf{Input:} {$N$,{$\{\plant{k}{k-1}\}_{k\in 1:N}$},$q_{0:N}$,$\{c_{k}^{\textup{x}}(\cdot)\}_{k\in1:N}$,$\{c_{k}^{\textup{u}}(\cdot)\}_{k\in1:N}$,$\{ \primitive{k}{k-1}{i}\}_{i \in \mathcal{P}}$}
\State \textbf{Output:} $\left\{\optimalpolicy{k}{k-1}\right\}_{1:N}$
\State $\costl{N+1}{N}\gets 0$
\For{$k = N:1$}
\State  $\bar{c}(\bv{X}_k,\bv{U}_k) \gets \costx{k} + \costu{k} + \costl{k+1}{k}$
\State $\optimalweights{k} \gets$ minimizer of the problem in \eqref{prob_stepk}
\State $\costl{k}{k-1} \gets $ optimal value of the problem in \eqref{prob_stepk}
\State $\optimalpolicy{k}{k-1}\gets\sum_{i\in\mathcal{P}}\optimalweight{k}{i}\primitive{k}{k-1}{i}$
\EndFor
\end{algorithmic}
\end{algorithm}

\subsection{Optimal Primitives Composition: Tackling Problem \ref{prob:finite-horizon}}\label{sec:properties}
In \Model, the free energy gating solves Problem \ref{prob:finite-horizon} and this is done via Algorithm \ref{alg_s2}, which outputs the optimal policy obtained by dynamically combining the primitives at each $k$. Given the inputs specified in Algorithm \ref{alg_s2}, it computes the optimal weights $\{\optimalweights{k}\}_{1:N}$ via the backward recursion in lines $4 - 9$. Within the recursion, $\optimalweights{k}$ is the minimizer of the following problem:
\begin{equation}\label{prob_stepk}
\begin{split}
    \min_{\weights{k} \in \Delta^{\np}} \ & \DKL{\jointxu{k}{k-1}}{\refjointxu{k}{k-1}} {+} \E_{\jointxu{k}{k-1}}\left[\bar{c}(\bv{X}_k,\bv{U}_k)\right]\\
    \mbox{s.t. } & \policy{k}{k-1} = \sum_{i\in \mathcal{P}}\weight{k}{i} \primitive{k}{k-1}{i},
\end{split}
\end{equation}
where $\bar{c}(\bv{X}_k,\bv{U}_k):=\costx{k}{+}\costu{k}{+} \costl{k+1}{k}$; $\costl{k+1}{k}$ is initialized at $0$ and iteratively built within the recursion in Algorithm \ref{alg_s2}. 
\begin{Remark}
\label{rem:cost_to_go}
The cost functional in \eqref{prob_stepk} is the variational free energy at time step $k$. In this problem, the expected loss contains the term $\costl{k+1}{k}$, which, as shown in Algorithm \ref{alg_s2}, 
is the optimal free energy value obtained by solving the problem at $k+1$. By embedding $\costl{k+1}{k}$ into the expected loss at time step $k$, the optimal policy $\optimalpolicy{k}{k-1}$ accounts for planning over the time horizon. For biomimetic (or greedy) actions, where $N=1$,  $\costl{k+1}{k}$ is equal to $0$.
\end{Remark}
Next, we characterize the optimality of the policy returned by Algorithm \ref{alg_s2} (Property \ref{prop:optimality}) and the convexity of the optimization problems solved at each $k$ (Property \ref{prop:convexity}).

\subsubsection{Properties of Algorithm \ref{alg_s2}}
We now analyze key properties of Algorithm~\ref{alg_s2}. Specifically, we establish the optimality of the policy returned by Algorithm \ref{alg_s2} and the convexity of the optimization problems that yield the optimal weights. The results leverage standard arguments from convex functions and control the space of densities. See Appendix~\ref{app:proof_optimality} for the self-contained proofs.

\begin{Property}[Optimality]\label{prop:optimality}
Algorithm \ref{alg_s2} returns the optimal solution of Problem \ref{prob:finite-horizon}.
\end{Property}
The next property characterizes the convexity of the problem in \eqref{prob_stepk}. The property can be proved by leveraging convexity of the KL divergence and the composition property of a convex function with an affine mapping. In Appendix~\ref{app:proof_convexity} we provide a proof by explicitly computing the Hessian of the cost function.
\begin{Property}[Convexity]\label{prop:convexity}
At each $k$, the problem in \eqref{prob_stepk} is convex in the decision variables $\weights{k}$. Moreover, if the primitives are linearly independent the cost is strictly convex.
\end{Property}
Remarkably, when the primitives are linearly independent, the problem features a strictly convex cost function. This implies that, if the optimal solution is in the relative interior of the feasibility domain, then \Model~policy outperforms the individual primitives, thus achieving {\em transcendence} in the spirit of \cite{zhang2024transcendence}.

\begin{Remark}[Transcendence]\label{rem:transcendence}
If the primitives are linearly independent and the optimal weight vector $\optimalweights{k} \in \Delta^{\np}$ has at least two strictly positive elements, then the resulting policy $\optimalpolicy{k}{k-1}$ strictly outperforms each primitive, i.e., $J_k(\optimalpolicy{k}{k-1}) < \ds \min_{i \in \mathcal{P}} J_k(\primitive{k}{k-1}{i})$, where $J_k(\cdot)$ denotes the cost in \eqref{prob_stepk}.
\end{Remark}

\section{Validation}\label{sec:validation}
\begin{figure}[tb]
\centering
    \subfigure[]{\includegraphics[width=0.45\columnwidth]{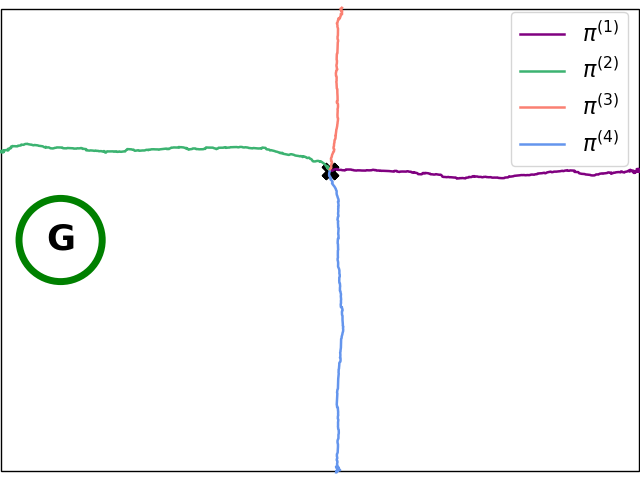}
    \label{robot_primitives}}
    \subfigure[]{\includegraphics[width=0.5\columnwidth]{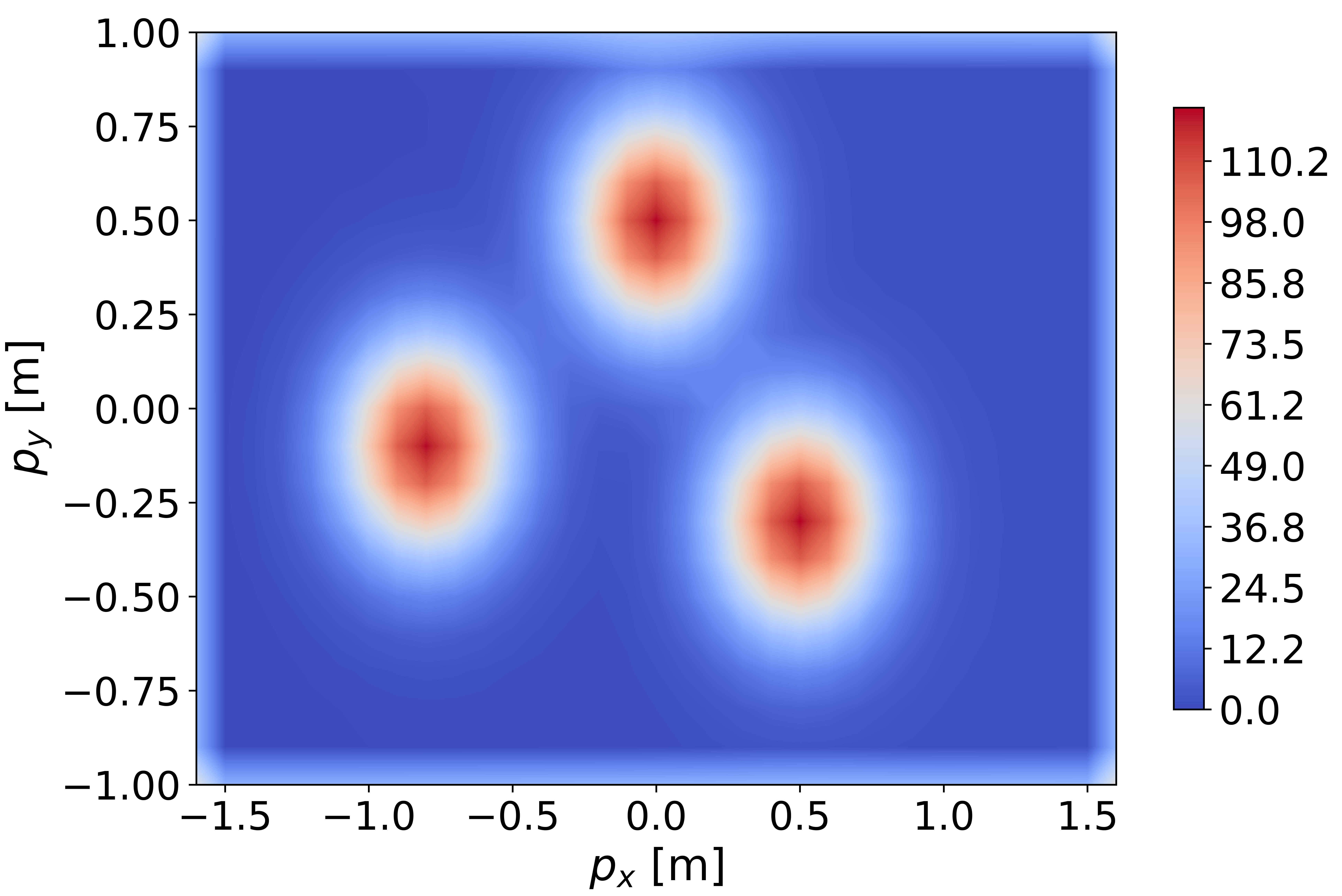}
    \label{robot_cost_map}} 
    \subfigure[]{\includegraphics[width=0.45\columnwidth]{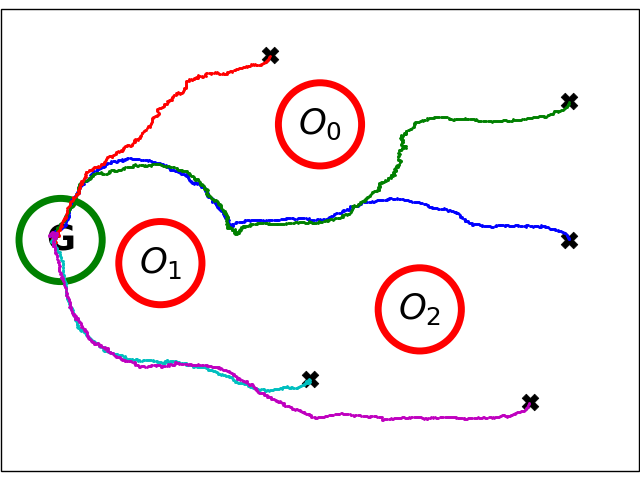}
    \label{robot_combination_obs}}
    \subfigure[]{\includegraphics[width=0.5\columnwidth]{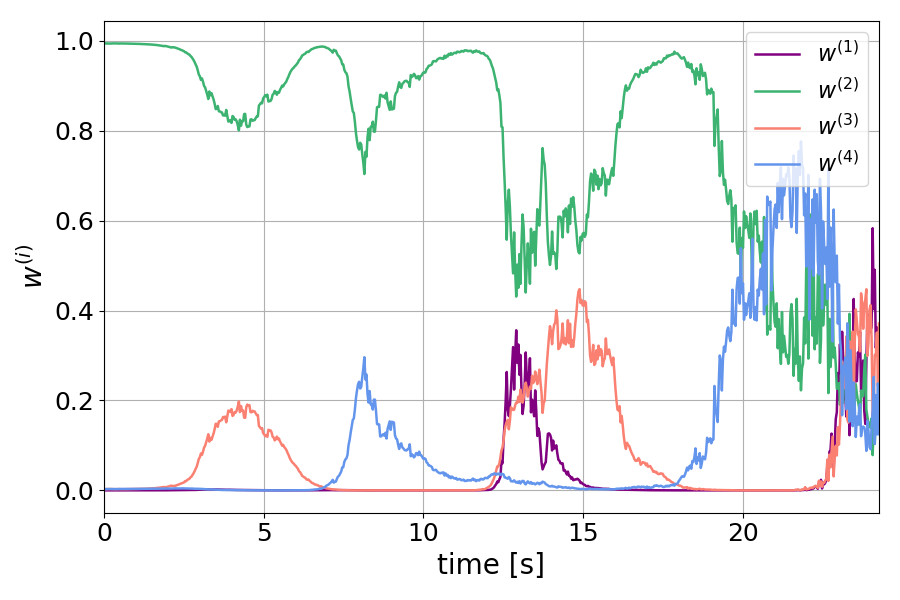}
    \label{robot_combination_weights}}
    \caption{(a) trajectories when the robot is controlled by the $4$ primitives described in Section \ref{sec:validation}; (b) heat map of the cost used in the experiments; (c) robot trajectories when controlled by \Model; (d) optimal primitives' weights corresponding to the experiment in blue in panel (c).
    {\Model~allows the robot to idle at the goal position despite lacking a stopping primitive.}
    }
\end{figure}
We use \Model~to control a unicycle robot so that it can navigate towards a goal point in an environment with obstacles. \Model~is validated through the Robotarium \cite{robotarium}: this platform offers both a hardware infrastructure and a high-fidelity simulator. Consequently, \Model~is validated via both simulations and real hardware experiments with rovers moving in a $\SI{3.2}{\meter} \times \SI{2}{\meter}$ work area. The code for replicating our experiment is available at \url{https://github.com/freegate-iwai/Free-Gate.git} together with a recording of our hardware experiment. The Robotarium provides built-in functions to map the unicycle dynamics onto a single integrator. Hence, it gives users the possibility to maneuver the robot considering an integrator dynamics. The dynamics we use is therefore $\bv{x}_k = \bv{x}_{k-1} + \bv{u}_k dt$, where $\bv{x}_k = [p_{x,k}, p_{y,k}]^{\top}$ is the state (robot's position, in \SI{}{\meter}) at time step $k$, $\bv{u}_k = [v_{x,k}, v_{y,k} ]^{\top}$ is the input vector (velocity, in \SI{}{\meter\per\second}) and $dt = \SI{0.033}{\second}$ is the Robotarium time step. We set $\bv{x}_k \in \mathcal{X}:=[-1.6, 1.6] \times [-1, 1]$, matching the work area, and $\bv{u}_k \in \mathcal{U} := [-0.2, 0.2] \times [-0.2,0.2]$, in accordance with the maximum speed physically allowed on the robot. For the application of our results, the state and input spaces are discretized into $33\times21$ and $7\times7$ equidistant bins, respectively. We emulate measurement noise for the robot position, which following, e.g., \cite{yoo2016mapless} we assume to be Gaussian. Hence, $\plant{k}{k-1} = \mathcal{N}(\bv{x}_{k-1} + \bv{u}_k dt, \boldsymbol{\Sigma}_{x})$, with covariance matrix $\boldsymbol{\Sigma}_{x}=0.008\bv{I}_2$ ($\bv{I}_2$ is the $2 \times 2$ identity). \Model~has access to $\np=4$ primitives, {which we obtained by discretizing Gaussians of the form $\primitive{k}{k-1}{i} = \mathcal{N}(\bar{\bv{u}}_k^{(i)}, \boldsymbol{\Sigma}_u)$. For these Gaussians, we set $\boldsymbol{\Sigma}_u= 0.005 \bv{I}_2$, while $\bar{\bv{u}}_k^{(i)}$ is a vector signal (e.g., coming from a simple proportional controller) proportional to the distance between the robot position and each of the boundaries of the Robotarium work area (see GitHub for details). These primitives only allow the robot to move along four cardinal directions (right, left, up, down -- see Fig. \ref{robot_primitives}) but cannot solve the navigation task individually.  Inspired by, e.g., \cite{CH_etal:24} and references therein, in the experiments we set the generative, time-series, model as: (i) $\refplant{k}{k-1}=\mathcal{N}(\bv{x}_{k-1} + \bv{u}_k dt, \tilde{\boldsymbol{\Sigma}}_x)$, with $\tilde{\boldsymbol{\Sigma}}_x = 0.002 \bv{I}_2$, and (ii) $\refpolicy{k}{k-1}= \mathcal{N}(\tilde{\bv{u}}_k,\tilde{\boldsymbol{\Sigma}}_{u})$, where $\tilde{\bv{u}}_k$ is provided by the Robotarium built-in position controller, not accounting for obstacles, and $\tilde{\boldsymbol{\Sigma}}_{u}=0.005\bv{I}_2$. For the control task, the goal position is $\subscr{\bv{x}}{d}=[-1.3, 0]^{\top}$ (Fig. \ref{robot_primitives}) and the robot needs to avoid obstacles positioned at $[0, 0.5]^{\top}$, $[-0.8, -0.1]^{\top}$, $[0.5, -0.3]^{\top}$ and captured via the following cost adapted from \cite{Garrabe2025}:
$c_{k}^{\textup{x}}(\bv{x}_k) = 150 \sum_{j=1}^3 g_j(\bv{x}_k) + 30 \sum_{j=1}^4 w_j(\bv{x}_k)$.
In this expression:
(i) $g_j(\cdot)$ are Gaussian-shaped obstacles penalties;
(ii) $w_j(\cdot)$ is the cost associated to the boundaries, with $w_j(\cdot)$ again defined as Gaussian (see the GitHub for the details). In Fig. \ref{robot_cost_map} the heat map of the cost $c_{k}^{\textup{x}}(\bv{x}_k)$ is shown. {In addition, to equip the agent with planning abilities without implementing the full backward recursion in Algorithm \ref{alg_s2}, we find a heuristically approximated cost-to-go (Remark \ref{rem:cost_to_go}) as the cost obtained by applying twice the same input (see GitHub for details). 
Note that Assumption \ref{assum_primitives}, Assumption \ref{asn:bounded_cost} {and Assumption \ref{assum_support}} are all satisfied (Remark \ref{rem:assumptions}). In the experiments, CVX \cite{cvxpy} is used to solve the optimization problem in Algorithm \ref{alg_s2} at each time step. To validate \Model, we run $5$ experiments starting from different initial positions of the robot. {Each simulation terminates when the robot remains within a $\SI{0.08}{\meter}$ radius of the goal for at least $\SI{2}{\second}$}. {We record the robot's trajectory in each experiment, and the results are shown in Fig. \ref{robot_combination_obs}: \Model~allows the robot to successfully reach the goal while avoiding obstacles. {Note that there is no dedicated {\em stop} primitive, yet \Model~still enables the robot to remain at the goal point once this is reached by properly combining the available primitives.}} This is also illustrated in Fig. \ref{robot_combination_weights}, where the time evolution of the optimal weights from one of the experiments is reported}. The simulation results are confirmed when \Model~is deployed on the Robotarium hardware. A recording of one hardware experiment is available on our GitHub. Namely, the recording shows that, when \Model~is used, the real robot successfully completes the navigation task. 

\section{Conclusions}\label{sec:conclusions}

We introduced \Model: a computational model for planning and control via policy composition. \Model~consists of a  set of functional units that return control primitives. The primitives are composed by a gating mechanism that minimizes the variational free energy. This leads to recast the problem of finding the optimal weights for the primitives as a finite-horizon optimal control problem. This problem is convex even when the agent cost is not and when the environment is nonlinear, stochastic and non-stationary. After turning the results into an algorithm, the effectiveness of \Model~is illustrated via both in-silico and hardware experiments on a robot navigation task. Based on the results presented here, our future work will be aimed at: (i) expanding \Model~to integrate, in the spirit of active inference, perception and learning; (ii) giving rigorous conditions establishing what makes for a good set of primitives; (iii) embedding in \Model~function approximators, e.g., deep neural networks, to estimate the cost-to-go and tackle continuous-control settings; (iv) embodying \Model~in applications where primitives are the output of modules are engineered to process multi-stream, heterogeneous, sensory information; (v) implementing the gating mechanism presented here via biologically plausible neural networks using, e.g., the tools inspired by \cite{centorrino_positive_2024}. 
\begin{credits}
\subsubsection{\ackname} 
GR wishes to thank K. Friston (University College London, UK) for the insightful discussions on the free energy minimization formulation. GR is also grateful to G. Pezzulo, D. Maisto (National Research Council, Italy), L. Da Costa (VERSES AI Research, USA), A. Paul (Monash University, Australia) for the conversations on the ideas of this paper.
\subsubsection{\discintname}
The authors have no competing interests to declare.
\end{credits}

\newpage
\appendix
\label{appendix}
\section{Proofs}
\paragraph{\bf Proof of Proposition \ref{prop:optimality}.}
\label{app:proof_optimality}
The derivations are inspired by \cite{DG_GR:22} and leverage the chain rule for the KL divergence. For notational convenience, in what follows we omit the subscript of the expectation operator when it applies to $\bar{c}(\bv{X}_N,\bv{U}_N)$, writing $\E\left[\bar{c}(\bv{X}_N,\bv{U}_N)\right]$ instead of $\E_{\jointxu{N}{N-1}}\left[\bar{c}(\bv{X}_N,\bv{U}_N)\right]$.
Specifically, following the chain rule, the cost in \eqref{eqn:finite-horizon} can be written as:
\begin{align*}
&\DKL{p_{0:N-1}}{q_{0:N-1}} + \sum_{k=1}^{N-1}\E_{\marginalx{k-1}}\left[\E_{\jointxu{k}{k-1}}\left[\costx{k}{+}\costu{k}\right]\right] \\
&{{+}\E_{\marginalx{N-1}}\left[\DKL{\jointxu{N}{N-1}}{\refjointxu{N}{N-1}}{+}\E\left[\bar{c}(\bv{X}_N,\bv{U}_N)\right]\right]},
\end{align*}
where $\bar{c}(\bv{X}_N,\bv{U}_N){=} \costx{N} {+} \costu{N} {+}\costl{N+1}{N}$, with $\costl{N+1}{N} {=} 0$.
Hence, Problem \ref{prob:finite-horizon} can be recast as the sum of the following two sub-problems:
\begin{equation*}
\begin{aligned}
    {\min_{\{\weights{k}\in \Delta^{\np}\}_{1:N-1}}} \!\!\!\! & {\DKL{p_{0:N-1}}{q_{0:N-1}}\!{+}\!\! \sum_{k=1}^{N-1}\!\E_{\marginalx{k-1}} \! \left[\E_{\jointxu{k}{k-1}} \! \left[\costx{k}\!{+}\costu{k}\right]\right]}\\
    \mbox{s.t. } &\policy{k}{k-1} = \sum_{i\in \mathcal{P}}\weight{k}{i}\primitive{k}{k-1}{i}, \ \ \ \forall k \in 1:N-1,
\end{aligned}
\end{equation*}
and
\begin{equation*}
\begin{split}
    \min_{\weights{N}\in \Delta^{\np}} & {\E_{\marginalx{N-1}}\left[\DKL{\jointxu{N}{N-1}}{\refjointxu{N}{N-1}} {+} \E\left[\bar{c}(\bv{X}_N,\bv{U}_N)\right]\right]}\\
   \mbox{s.t. } & \policy{N}{N-1} = \sum_{i\in \mathcal{P}}\weight{N}{i}\primitive{N}{N-1}{i},
\end{split}
\end{equation*}
In the above problem, since the expectation in the cost functional does not depend on $\bv{w}_N$, the optimal cost is $\E_{\marginalx{N-1}}[\costl{N}{N-1}]$, where $\costl{N}{N-1}$ is obtained by solving:
\begin{equation*}
\begin{split}
\min_{\weights{N}\in \Delta^{\np}} & \DKL{\jointxu{N}{N-1}}{\refjointxu{N}{N-1}} + \E\left[\bar{c}(\bv{X}_N,\bv{U}_N)\right]\\
  \mbox{s.t. } & \policy{N}{N-1} = \sum_{i\in \mathcal{P}}\weight{N}{i}\primitive{N}{N-1}{i},
\end{split}
\end{equation*}
This is the optimization problem solved at $k=N$ in Algorithm \ref{alg_s2}.\\
Note that $\E_{\marginalx{N-1}}[\costl{N}{N-1}]$ is bounded (Assumption \ref{asn:bounded_cost}) and can be conveniently written as $\E_{\marginalx{N-2}}\bigg[\E_{\jointxu{N-1}{N-2}}\big[\costl{N}{N-1}\big]\bigg]$, so that the cost in \eqref{eqn:finite-horizon} becomes $\DKL{p_{0:N-1}}{q_{0:N-1}} {+} \ds {\sum_{k=1}^{N-1}}\E_{\marginalx{k-1}}\bigg[\E_{\jointxu{k}{k-1}}\big[\costx{k}+\costu{k}\big]\bigg]
+ \E_{\marginalx{N-2}}\bigg[\E_{\jointxu{N-1}{N-2}}\big[\costl{N}{N-1}\big]\bigg]$. 
Hence, iterating the arguments above, we obtain that the optimal solution of Problem \ref{prob:finite-horizon} can be found by solving, at each $k$, the optimization problem given in Algorithm \ref{alg_s2} with $\costl{N+1}{N}$ equal to $0$ and, $\forall k\in 1:N-1$, $\costl{k+1}{k}$ being the optimal value of the optimization problem at $k+1$. This yields the desired conclusions. \qed

\paragraph{\bf Proof of Property~\ref{prop:convexity}.} 
\label{app:proof_convexity}
Following \cite{garrabe2023optimal}, we note that the constraint set is convex and therefore we only need to prove convexity of the cost. We do it by showing that the Hessian, $\mathbf{H}(\bv{x}_{k-1})$, with entries $h_{ij}$ for $i,j \in \mathcal{P}$, is positive semi-definite. 
For compactness, we introduce the following shorthand notation: $\pi_{\bv{u},k}^{(i)}:=\primitive{k}{k-1}{i}$, $p_{\bv{x},k}:=\plant{k}{k-1}$, $q_{\bv{x,u},k}:=\refjointxu{k}{k-1}$.
Embedding the constraint $\policy{k}{k-1} = \sum_{i\in \mathcal{P}}\weight{k}{i} \primitive{k}{k-1}{i}$ directly into the cost functional of problem \eqref{prob_stepk}, $
l_k(\bv{x}_{k-1}):=\DKL{\jointxu{k}{k-1}}{\refjointxu{k}{k-1}} + \E_{\jointxu{k}{k-1}}\left[\bar c(\bv{X}_{k},\bv{U}_k)\right]$,
the first derivative of $l_k(\bv{x}_{k-1})$
with respect to the decision variable $\weight{k}{i}$ is:
\begin{equation*}
\begin{aligned}
     \frac{\partial l_k(\bv{x}_{k-1})}{\partial \weight{k}{i}} &{=} \frac{\partial}{\partial \weight{k}{i}} \Biggl( {\int} {\int} \sum_{j\in\mathcal{P}}\weight{k}{j}\pi_{\bv{u},k}^{(j)}p_{\bv{x},k}
     \ln \left(\frac{ \sum_{j\in\mathcal{P}}\weight{k}{j}\pi_{\bv{u},k}^{(j)}p_{\bv{x},k}}{q_{\bv{x,u},k}} \right) d\mathbf{x}_{k}d \mathbf{u}_{k}\\
     &\quad \quad \quad \quad +  {\int} {\int} \sum_{j\in\mathcal{P}}\weight{k}{j}\pi_{\bv{u},k}^{(j)}p_{\bv{x},k}\bar{c}(\bv{x}_{k},\bv{u}_k) d\mathbf{x}_{k}d \mathbf{u}_{k}\Biggr)\\
    &= \int \int \pi_{\bv{u},k}^{(i)} p_{\bv{x},k}\left(\ln\left(\frac{\sum_{j\in\mathcal{P}}\weight{k}{j}\pi_{\bv{u},k}^{(j)}p_{\bv{x},k}}{q_{\bv{x,u},k}}\right) +\bar c(\bv{x}_{k},\bv{u}_k)+1\right)d\mathbf{x}_{k}d \mathbf{u}_{k},
\end{aligned}
\end{equation*}
where the integrals, well defined by Assumption \ref{assum_primitives} and \ref{assum_support}, are respectively over
$\sX$ and $\sU$. Next, we compute the second order derivative: 
\begin{equation*}
    \begin{aligned}
     h_{ij}(\bv{x}_{k-1})&:=\frac{\partial^2 l_k(\bv{x}_{k-1})}{\partial \weight{k}{i}\partial \weight{k}{j}}
     = \int\int  \frac{\pi_{\bv{u},k}^{(i)}\pi_{\bv{u},k}^{(j)}p_{\bv{x},k}}{\sum_{t\in\mathcal{P}}\weight{k}{t}\pi_{\bv{u},k}^{(t)}}d\bv{x}_{k} d\bv{u}_{k}, \ \ \ \forall i,j \in \mathcal{P}.
    \end{aligned}
\end{equation*}
Pick now a non-zero vector $\mathbf{v} \in \mathbb{R}^{\np}$. We get:
\begin{equation}
\label{eq:vHv}
    \mathbf{v}^\top\mathbf{H}(\bv{x}_{k-1})\mathbf{v} =\sum_{i,j}v_i v_j h_{ij} (\bv{x}_{k-1})
    = \int\int \mathbf{v}^\top\mathbf{A}(\bv{x}_k, \bv{u}_k, \bv{x}_{k-1})\mathbf{v}\ d\bv{x}_k d\bv{u}_k,
\end{equation}
where $\mathbf{A}(\bv{x}_k, \bv{u}_k, \bv{x}_{k-1}) := \frac{\ds p_{\bv{x},k}}{\ds \sum_{t\in\mathcal{P}}\weight{k}{t}\pi_{\bv{u},k}^{(t)}} \boldsymbol{\pi}_{\bv{u},k} \boldsymbol{\pi}_{\bv{u},k}^\top$, with  $\boldsymbol{\pi}_{\bv{u},k} := [\pi^{(1)}_{\bv{u},k}, \ldots, \pi^{(\np)}_{\bv{u},k}]^\top$.
Convexity then follows by noticing that $\mathbf{A}(\bv{x}_k, \bv{u}_k, \bv{x}_{k-1})$, is the outer product of vector $\boldsymbol{\pi}_{\bv{u},k}$ by itself rescaled by a non-negative fraction. Therefore, $\mathbf{A}(\bv{x}_k, \bv{u}_k, \bv{x}_{k-1})$ is positive semi-definite for all $\bv{x}_k, \bv{u}_k, \bv{x}_{k-1}$.

To establish strict convexity, we notice that the quantity in \eqref{eq:vHv} can be rewritten as:
$$
\mathbf{v}^\top\mathbf{H}(\bv{x}_{k-1})\mathbf{v}= \int \int \frac{p_{\bv{x},k}}{\sum_t \weight{k}{t} \pi^{(t)}_{\bv{u},k}} (\bv{v}^\top \boldsymbol{\pi}_{\bv{u},k})^2 \, d\bv{x}_k d\bv{u}_k,
$$
which is zero if and only if $\mathbf{v}^\top \boldsymbol{\pi}_{\bv{u},k}= 0$ almost everywhere in the support of $p_{\bv{x},k}$. If the primitives $\pi^{(i)}_{\bv{u},k}$, which are positive over $\sU$ (by Assumption \ref{assum_support}), are linearly independent, there is no non-zero vector $\mathbf{v}$ such that $\mathbf{v}^\top \boldsymbol{\pi}_{\bv{u},k} = 0$~\cite{yu2018efficient}. Hence, in this case $\mathbf{H}(\bv{x}_{k-1})$ is positive definite and thus the cost functional $l_k(\bv{x}_{k-1})$ is strictly convex in $\weights{k}$. \qed

%
%
%
\newpage
\bibliographystyle{splncs04}

\bibliography{biblio}

\end{document}